\documentclass[12pt]{article}
\usepackage{amsmath}

\usepackage{t1enc}
\usepackage[latin1]{inputenc}
\usepackage[english]{babel}
\usepackage{amsmath,amsthm}
\usepackage{amsfonts}
\usepackage{latexsym}
\usepackage[dvips]{graphicx}
\usepackage{graphicx}
\usepackage[natural]{xcolor}
\newtheorem{theorem}{Theorem}
\newtheorem{lemma}[theorem]{Lemma}
\newtheorem{proposition}[theorem]{Proposition}

\begin{document}

\date{}
\author{Liangquan ZHANG \ and \ \ Yufeng SHI \thanks{
Corresponding author. E-mail: yfshi@sdu.edu.cn. Supported by National
Natural Science Foundation of China Grant 10771122, Natural Science
Foundation of Shandong Province of China Grant Y2006A08 and National Basic
Research Program of China (973 Program, No. 2007CB814900).} \\
School of Mathematics, Shandong University\\
Jinan 250100, People's Republic of China}
\title{Comparison Theorems of Infinite Horizon Forward-Backward Stochastic
Differential Equations }
\maketitle

\begin{abstract}
By the methods of probability and duality technique, we give some comparison
theorems for the solutions of infinite horizon forward-backwad stochastic
differential equations.
\end{abstract}


\section{Introduction}

In this paper, we study the comparison theorems of the following infinite
horizon fully coupled forward-backward stochastic differential equations
(FBSDEs in short)
\begin{equation}
\left\{
\begin{array}{l}
dX_t=b\left( t,X_t,Y_t,Z_t\right) dt+\sigma \left( t,X_t,Y_t,Z_t\right) dB_t,
\\
-dY_t=f\left( t,X_t,Y_t,Z_t\right) dt-Z_tdB_t, \\
X_0=x,\text{ \quad }Y_\infty =\Phi \left( X_\infty \right) ,
\end{array}
\right.  \tag{1.1}
\end{equation}
where $\left\{ B_t\right\} _{t\geq 0}$ be a standard $d$-dimensional
Brownian motion, defined on a probability space $\left( \Omega ,{\cal F}%
,P\right) $. Let $\left\{ {\cal F}_t\right\} _{t\geq 0}$ be the completion
of the natural filtration generated by the Brownian motion $\left\{
B_t\right\} _{t\geq 0},$ where ${\cal F}_0$ contains all $P$-null sets of $%
{\cal F}$, and $\mathcal{F}_\infty =\bigvee\limits_{t\geq
0}\mathcal{F}_t$
 We introduce the following spaces:
\begin{eqnarray*}
\mathcal{S}^2 &=&\left\{ \upsilon _t,\text{ }0\leq t\leq \infty ;\text{ }%
\upsilon _t\text{ is an }\mathcal{F}_t\text{-adapted process s.t. }\mathbf{E}%
\left[ \sup\limits_{0\leq t\leq \infty }\left| \upsilon _t\right|
^2\right]
<\infty .\right\} , \\
\mathcal{H}^2 &=&\left\{ \upsilon _t,\text{ }0\leq t\leq \infty ;\text{ }%
\upsilon _t\text{ is an }\mathcal{F}_t\text{-adapted process s.t. }\mathbf{E}%
\left[ \int_0^\infty \left| \upsilon _t\right| ^2dt\right] <\infty
.\right\}
, \\
L^2 &=&\left\{ \xi ;\text{ }\xi \text{ is an }\mathcal{F}_\infty \text{%
-measurable random variable s.t. }\mathbf{E}\left| \xi \right|
^2<\infty .\right\} ,
\end{eqnarray*}

and denote that ${\cal B}^3{\cal =S}^2{\cal \times S}^2{\cal \times H}^2.$

It is well known that the Hamiltonian system arised in studying maximum
principle of stochastic optimal control problems is a kind of fully coupled
FBSDEs (see [14]). In mathematical finance, fully coupled FBSDEs can be
encountered when one studied the problems of hedging options involved in a
large investor in financial market (see [2, 6]). On the other hand, fully
coupled FBSDEs can provide probabilistic interpretations for the solutions
of a class of quasilinear parabolic and elliptic PDEs (see [13, 16, 19]) and
have been investigated deeply. They were studied first by Antonelli [1] and
Ma, Protter and Yong [9], Cvitanic and Ma [2], Duffie, Ma and Yong [5].
Consequently, FBSDEs were developed in Hu and Peng [7], Pardoux and Tang
[13], Peng and Wu [14] and Yong [21].

In addition, infinite horizon FBSDEs play a role of studying the Black
consol rate conjecture (see [5]). Peng and Shi [15] firstly proved the
existence and uniqueness of solutions of infinite horizon FBSDEs, under some
monotone assumptions, while the solutions are in a kind of infinite time
integrable space and the infinite terminal values of the solutions must be
decay. Recently Shi and Zhao [16] extended the results in [15] to a larger
space where the solutions must not be decay in infinite horizon. Later Wu
[18] gave an existence and uniqueness result for FBSDEs with stopping times
under some similar monotone assumptions, but different Lipschitz conditions.
When the stopping times take infinite value, the solutions in [18] are in
the space ${\cal B}^3$ above and not necessarily decay. We will present the
existence and uniqueness result of [18] in Section 2.

The comparison theorems of FBSDEs play an important role of the applications
of FBSDEs. So many literatures discussed the comparison results for FBSDEs.
By virtue of partial differential equation method, Cvitanic and Ma [2]
obtained a comparison theorem of FBSDEs. They required the forward equation
to be non-degenerate and the coefficients not to be random, i.e. in
Markovian cases. Obviously their results can not treat the non-Markovian
cases that the coefficients themselves are randomly distributed, which are
often the fact in practice, e.g. in incomplete finance markets. Later for
non-Markovian finite horizon FBSDEs, Wu [17] proved a comparison result at
initial time, where the dimension of the forward system $X$ is $n>1$ and the
one of the backward system $Y$ is $m=1.$ Peng and Shi [15] gave a comparison
theorem for the solutions $X_t$ and $Y_t$ of infinite horizon FBSDEs for any
$t\in R^{+}$ with $n=m=1$. Recently Wu and Xu [20] obtained some comparison
results for FBSDEs on finite horizon for many situations, such as $m=1,$ $%
n>1 $ and $m>1,$ $n=1.$ Moreover, some results hold for the process $Y_t$ on
the whole interval of $t\in \left[ 0,T\right] $. However there is a little
gap in [20]. For overcoming this gap, we import a condition (i.e. (H4), for
details see the proof of Theorem 2 in this paper). We will extend all the
results in [17, 20, 22] to infinite horizon FBSDEs in many different
situations including multi-dimensional cases. In addition, we also extend
the comparison theorem for one-dimension infinite horizon backward
stochastic differential equations (BSDEs in short) given in [8] to
multi-dimensional case. These comparison theorems are significant in the
theories of stochastic optimal control, differential games, PDEs and
mathematical finance and so on.

This paper is organized as follows: we give the preliminaries and
assumptions in Section 2, all kinds of comparison theorems are proved in
Section 3 by virtue of probability method.

\section{Preliminaries: the existence and uniqueness to FBSDEs}

We consider FBSDE (1.1) under the following setting: for any $t\geq 0$, $%
\left( X,Y,Z\right) \in R^n\times R^m\times R^{m\times d},$%
\begin{eqnarray*}
b\left( t,X,Y,Z\right) &:&\Omega \times \left[ 0,\infty \right] \times {\bf R%
}^n{\bf \times R}^m{\bf \times R}^{m\times d}{\bf \rightarrow R}^n, \\
\sigma \left( t,X,Y,Z\right) &:&\Omega \times \left[ 0,\infty \right] \times
{\bf R}^n{\bf \times R}^m{\bf \times R}^{m\times d}{\bf \rightarrow R}%
^{n\times d}, \\
f\left( t,X,Y,Z\right) &:&\Omega \times \left[ 0,\infty \right] \times {\bf R%
}^n{\bf \times R}^m{\bf \times R}^{m\times d}{\bf \rightarrow R}^m, \\
\Phi \left( X\right) &:&\Omega \times {\bf R}^n{\bf \rightarrow R}^m.
\end{eqnarray*}

We assume

\begin{enumerate}
\item[(H1)]  $\forall \left( X,Y,Z\right) \in {\bf R}^{n+m+m\times d},$ $%
\Phi \left( X\right) \in {\bf L}^2,$ $b,$ $\sigma ,$ and $f$ are
progressively measurable and
\[
\mathbf{E}\left( \int_0^\infty \left| b(s,0,0,0)\right| ds\right) ^2+\mathbf{%
E}\left( \int_0^\infty \left| f(s,0,0,0)\right| ds\right) ^2+\mathbf{E}%
\int_0^\infty \left| \sigma \left( s,0,0,0\right) \right|
^2ds<\infty ;
\]

\item[(H2)]  There exists a positive deterministic bounded function $u_1(t),$
such that $\forall $ $t\geq 0$, $\forall \left( X_i,Y_i,Z_i\right) \in {\bf R%
}^{n+m+m\times d},$ ($i=1,2$),
\[
\left| l\left( t,X_1,Y_1,Z_1\right) -l\left( t,X_2,Y_2,Z_2\right) \right|
\leq u_1\left( t\right) \left( \left| X_1-X_2\right| +\left| Y_1-Y_2\right|
+\left| Z_1-Z_2\right| \right) ,
\]
$l=b,$ $\sigma ,$ $f,$ respectively, and $\int_0^\infty u_1\left(
s\right) ds<\infty ,$ $\int_0^\infty u_1^2\left( s\right) ds<\infty
$
 there exists a positive
constant $C$ such that
\[
\left| \Phi \left( X_1\right) -\Phi \left( X_2\right) \right| \leq C\left|
X_1-X_2\right| ,\quad u_1\leq C.
\]
\end{enumerate}

Given an $m\times n$ full rank matrix $G$ and let us introduce the following
notations
\[
u=\left(
\begin{array}{c}
X \\
Y \\
Z
\end{array}
\right) ,A\left( t,u\right) =\left(
\begin{array}{c}
-G^{*}f \\
Gb \\
G\sigma
\end{array}
\right) \left( t,u\right) ,
\]
where $G\sigma =\left( G\sigma _1\ldots G\sigma _d\right) $. We use the
usual inner product and Euclidean norm in ${\bf R}^n{\bf ,R}^m{\bf ,R}%
^{m\times d}.$ All the equalities and inequalities mentioned in this
paper are in sense of $dt\times dp$ almost surely on $\left[
0,\infty \right] \times \Omega .$ In addition, throughout this
paper, when we indicate two vectors $\nu ^1,$ $\nu ^2$ satisfying
$\nu ^1\geq \nu ^2$, it means that we have $\nu ^{1,j}\geq \nu
^{2,j}$ for each of their components.
\begin{enumerate}
\item[H3]  $\forall u=\left( X,Y,Z\right) ,$ $\bar u=\left( \bar X,\bar
Y,\bar Z\right) \in \mathbf{R}^{n+m+m\times d},$ $\hat X=X-\bar X,$
$\hat
Y=Y-\bar Y,$ $\hat Z=Z-\bar Z,$%
\[
\left\{
\begin{array}{l}
\left\langle A\left( t,u\right) -A\left( t,\bar u\right) ,u-\bar
u\right\rangle \leq -\beta _1u_1\left( t\right) \left| G\hat
X\right| ^2-\beta _2u_1\left( t\right) \left( \left| G^{*}\hat
Y\right| ^2+\left|
G^{*}\hat Z\right| ^2\right) , \\
\left\langle \Phi \left( X\right) -\Phi \left( \bar X\right)
,G\left( X-\bar X\right) \right\rangle \geq \mu \left| G\hat
X\right| ^2,
\end{array}
\right.
\]
or
\[
\left\{
\begin{array}{l}
\left\langle A\left( t,u\right) -A\left( t,\bar u\right) ,u-\bar
u\right\rangle \geq \beta _1u_1\left( t\right) \left| G\hat X\right|
^2+\beta _2u_1\left( t\right) \left( \left| G^{*}\hat Y\right|
^2+\left|
G^{*}\hat Z\right| ^2\right) , \\
\left\langle \Phi \left( X\right) -\Phi \left( \bar X\right)
,G\left( X-\bar X\right) \right\rangle \leq -\mu \left| G\hat
X\right| ^2,
\end{array}
\right.
\]
\end{enumerate}
where $\beta _1,$ $\beta _2$ and $\mu ,$ are given nonnegative constant with
$\beta _1+\beta _2>0,$ $\mu +\beta _2>0.$ Moreover, we have $\beta _1>0,$ $%
\mu >0$ (resp., $\beta _2>0$) when $m>n$ (resp., $n>m$). Particularly, when $%
m=n$, we set $G=I_n.$ And we denote the transpose of matrix by the
notation *.

\begin{enumerate}
\item[H4]  $f$ is strictly monotone in $X$.
\end{enumerate}
\begin{proposition}
Assume (H1)-(H3) hold, then FBSDE (1.1) has a unique solution $\left(
X\left( \cdot \right) ,Y\left( \cdot \right) ,Z\left( \cdot \right) \right)
\in {\cal B}^3.$
\end{proposition}

The proof can be seen in [18].

\bigskip
For convenience of notations, hereafter we set $d=1$ in this paper. In order
to prove multi-dimensional comparison theorems for FBSDE (1.1), we need the
following so-called quasi-monotonicity conditions.

\begin{enumerate}
\item[(A1)]  For $\forall t\geq 0,$ $\forall Y\in {\bf R}^m,$ $Z{\bf \in R}%
^{m\times d},$ $b\left( t,\cdot ,Y,Z\right) :\Omega \times {\bf R}^n{\bf %
\rightarrow R}^n$ satisfies that $\forall X,$ $\bar X\in {\bf R}^n$, if $%
X_j\leq \bar X_j$ ($1\leq j\leq n$), while $X_s=\bar X_s$ ($s\neq j$), then $%
b^s\left( t,X,Y,Z\right) \leq b^s\left( t,\bar X,Y,Z\right) $. For $\forall
Y\in {\bf R}^m,$ $Z{\bf \in R}^m$, $\sigma \left( t,\cdot ,Y,Z\right)
:\Omega \times {\bf R}^{+}{\bf \times R}^n{\bf \rightarrow R^n},$ satisfies
that $\sigma ^j\left( t,x\right) =\sigma ^j\left( t,0,0,\ldots ,x_j,\ldots
,0,Y,Z\right) .$

\item[(A2)]  For $\forall t\geq 0,$ $\forall X\in {\bf R}^n$, $f_s^i\left(
t,X,\cdot ,\cdot \right) :\Omega \times {\bf R}^n{\bf \times R}^n{\bf %
\rightarrow R}^n$ ($i$=$1,2$) satisfies that $\forall Y,$ $\bar Y\in {\bf R}%
^m,$ $\forall Z,$ $\bar Z\in {\bf R}^m,$ if $Y_s=\bar Y_s,$ $Z_s=\bar Z_s,$ (%
$s\neq j$), $Y_j\geq \bar Y_j,$ $(1\leq j\leq n)$, then $f_s^1\left(
t,X,Y,Z\right) \geq f_s^2\left( t,X,\bar Y,\bar Z\right) $.
\end{enumerate}

\section{Comparison theorem}

\subsection{1-dimensional cases of FBSDEs}

We first give the following result for $n\geq 1,$ $m\geq 1,$ which plays a
significant role in our following comparison results.

\begin{theorem}
Let $\left( X_{\cdot }^1,Y_{\cdot }^1,Z_{\cdot }^1\right) $ and $\left(
X_{\cdot }^2,Y_{\cdot }^2,Z_{\cdot }^2\right) $ be respectively the
solutions of (1.1) corresponding to $X_0^1=x_1\in {\bf R}^n$ and $%
X_0^2=x_2\in {\bf R}^n.$ Set
\begin{eqnarray*}
\hat U &=&U^1-U^2=\left( X^1,Y^1,Z^1\right) -\left( X^2,Y^2,Z^2\right) \\
&=&\left( X^1-X^2,Y^1-Y^2,Z^1-Z^2\right) =\left( \hat X,\hat Y,\hat Z\right)
.
\end{eqnarray*}
If (H1)-(H4) hold, then $\left\langle G\hat X_t,\hat Y_t\right\rangle \geq 0$%
, $\forall t\in R^{+}.$ Moreover, $\left( \hat X_t,\hat Y_t,\hat Z_t\right)
{\bf I}_{\left[ \tau ,\infty \right] }\left( t\right) \equiv 0$, where $\tau
$ is an ${\cal F}_t$-stopping time defined by
\[
\tau \doteq \inf \left\{ t\geq 0;\text{ }\left\langle G\hat X_t,\hat
Y_t\right\rangle =0\right\} .
\]
\end{theorem}

\noindent \textbf{Proof.}\quad \noindent For $\forall t\in {\bf
R}^{+}$, $\forall T\in [t,\infty ),$ applying It\^o's formula to
$\left\langle G\hat X,\hat Y\right\rangle $ on $[t,T],$ it follows
that
\begin{eqnarray*}
&&\mathbf{E}^{\mathcal{F}_t}\left\langle G\hat X_T,\hat
Y_T\right\rangle
-\left\langle G\hat X_t,\hat Y_t\right\rangle  \\
&=&\mathbf{E}^{\mathcal{F}_t}\int_t^T\langle A\left( s,U_s^1\right)
-A\left(
s,U_s^2\right) ,\hat U_s\rangle ds \\
&\leq &-\beta _1\mathbf{E}^{\mathcal{F}_t}\int_t^Tu_1(s)\left| G\hat
X_s\right| ^2ds-\beta _2\mathbf{E}^{\mathcal{F}_t}\int_t^Tu_1\left(
s\right) \left( \left| G^{*}\hat Y_s\right| ^2+\left| G^{*}\hat
Z_s\right| ^2\right) ds.
\end{eqnarray*}

Letting $T\rightarrow \infty ,$ and by (H3) it follows that
\begin{eqnarray*}
\left\langle G\hat X_t,\hat Y_t\right\rangle  &\geq &\mu \mathbf{E}^{%
\mathcal{F}_t}\left| G\hat X_\infty \right| ^2+\beta _1\mathbf{E}^{\mathcal{F%
}_t}\int_t^\infty u_1\left( s\right) \left| G\hat X_s\right| ^2ds \\
&&+\beta _2\mathbf{E}^{\mathcal{F}_t}\int_t^\infty u_1\left(
s\right) \left(
\left| G^{*}\hat Y_s\right| ^2+\left| G^{*}\hat Z_s\right| ^2\right) ds \\
&\geq &0.
\end{eqnarray*}

So we have $\left\langle G\hat X_t,\hat Y_t\right\rangle \geq 0,$ $\forall
t\geq 0.$ It follows that, for the above given stopping time $\tau $%

\begin{eqnarray*}
&&\mathbf{E}\left( \langle G\hat X_\infty ,\hat Y_\infty \rangle
-\langle
G\hat X_\tau ,\hat Y_\tau \rangle \right)  \\
&=&\mathbf{E}\int_\tau ^\infty \langle A\left( s,U_s^1\right)
-A\left(
s,U_s^2\right) ,\hat U_s\rangle ds \\
&\leq &-\beta _1\mathbf{E}\int_\tau ^\infty u_1\left( s\right)
\left| G\hat X_s\right| ^2ds-\beta _2\mathbf{E}\int_\tau ^\infty
u_1\left( s\right) \left( \left| G^{*}\hat Y_s\right| ^2+\left|
G^{*}\hat Z_s\right| ^2\right) ds,
\end{eqnarray*}

which means that
\[
\beta _1\mathbf{E}\int_\tau ^\infty u_1\left( s\right) \left| G\hat
X_s\right| ^2ds+\beta _2\mathbf{E}\int_\tau ^\infty u_1\left(
s\right) \left( \left| G^{*}\hat Y_s\right| ^2+\left| G^{*}\hat
Z_s\right| ^2\right) ds\equiv 0.
\]
When $m>n$, we have $\beta _1>0,$ so $\hat X_t{\bf I}_{\left[ \tau
,\infty \right] }\left( t\right) \equiv 0.$ From Theorem 2.1 in
[18], we get $\hat
Y_t{\bf I}_{\left[ \tau ,\infty \right] }\left( t\right) \equiv 0,$ $\hat Z_t%
{\bf I}_{\left[ \tau ,\infty \right] }\left( t\right) \equiv 0.$

When $m<n$, we consider the following FBSDEs
\[
\left\{
\begin{array}{l}
d\hat X_t=\left[ b_t^1\hat X_t+b_t^2\hat Y_t+b_t^3\hat Z_t\right] dt+\left[
\sigma _t^1\hat X_t+\sigma _t^2\hat Y_t+\sigma _t^3\hat Z_t\right] dB_t, \\
-d\hat Y_t=\left[ f_t^1\hat X_t+f_t^2\hat Y_t+f_t^3\hat Z_t\right] dt-\hat
Z_tdB_t, \\
\hat X_0=X_0^1-X_0^2,\quad \quad \hat Y_\infty =\Phi \left( X_\infty
^1\right) -\Phi \left( X_\infty ^2\right) ,
\end{array}
\right.
\]
where
\begin{eqnarray*}
l_t^{1j} &=&\left\{
\begin{array}{l}
\frac{l\left( t,X_t^{21},\ldots ,X_t^{1j},\ldots
,X_t^{1n},Y_t^1,Z_t^1\right) -l\left( t,X_t^{21},\ldots ,X_t^{2j},\ldots
,X_t^{1n},Y_t^1,Z_t^1\right) }{X_t^{1j}-X_t^{2j}},\text{ \ \ \ if }%
X_t^{1j}\neq X_t^{2j}, \\
0,\text{\qquad otherwise,}
\end{array}
\right. \\
l_t^{2j} &=&\left\{
\begin{array}{l}
\frac{l\left( t,X_t^2,Y_t^{21},\ldots ,Y_t^{1j},\ldots
,Y_t^{1m},Z_t^1\right) -l\left( t,X_t^2,Y_t^{21},\ldots ,Y_t^{2j},\ldots
,Y_t^{1m},Z_t^1\right) }{Y_t^{1j}-Y_t^{2j}},\text{\quad if }Y_t^{1j}\neq
Y_t^{2j}, \\
0,\qquad \text{otherwise,}
\end{array}
\right. \\
l_t^{3j} &=&\left\{
\begin{array}{l}
\frac{l\left( t,X_t^2,Y_t^2,Z_t^{21},\ldots ,Z_t^{1j},\ldots
,Z_t^{1m}\right) -l\left( t,X_t^2,Y_t^2,Z_t^{21},\ldots ,Z_t^{2j},\ldots
,Z_t^{1m}\right) }{Z_t^{1j}-Z_t^{2j}},\text{ \ \ \ if }Z_t^{1j}\neq Z_t^{2j},
\\
0,\text{\qquad otherwise.}
\end{array}
\right.
\end{eqnarray*}
Here $l=b,$ $\sigma ,$ $f,$ respectively. And note that
\[
\hat X_t=X_t^1-X_t^2=\left(
\begin{array}{c}
X_t^{11}-X_t^{21} \\
X_t^{12}-X_t^{22} \\
\vdots \\
X_t^{1n}-X_t^{2n}
\end{array}
\right) .
\]
$\hat Y_t,$ $\hat Z_t$ also have similar notation. Therefore, $b_t^1$ and $%
\sigma _t^1$ are $n\times n$ matrices. $b_t^2,$ $\sigma _t^2,$ $b_t^3$ and $%
\sigma _t^3$ are $n\times m$ matrices. $f_t^1$ is an $m\times n$ matrix. $%
f_t^2$ and $f_t^3$ are $m\times m$ matrices. When $m<n$, we have $\beta
_2>0, $ then $\hat Y_t{\bf I}_{\left[ \tau ,\infty \right] }\left( t\right)
\equiv 0,$ $\hat Z_t{\bf I}_{\left[ \tau ,\infty \right] }\left( t\right)
\equiv 0,$ Then we have
\[
\hat f_t^1\hat X_t{\bf I}_{\left[ \tau ,\infty \right] }\left( t\right)
\equiv 0,
\]
From (H4) it follows that $\hat X_\tau =0,$ a.s.. Using the uniqueness of
solutions of stochastic differential equations, it is obvious that $\hat X_t%
{\bf I}_{\left[ \tau ,\infty \right] }\left( t\right) \equiv 0.$ In all
\[
\hat X_t{\bf I}_{\left[ \tau ,\infty \right] }\left( t\right) \equiv \hat Y_t%
{\bf I}_{\left[ \tau ,\infty \right] }\left( t\right) \equiv \hat Z_t{\bf I}%
_{\left[ \tau ,\infty \right] }\left( t\right) \equiv 0.
\]
The proof is complete.\quad$\Box$

Now we first consider a simple case, that is, $n=m=1.$

\begin{theorem}
Assume that $n=m=1$ and $\left( b,\sigma ,f,\Phi \right) $ satisfy
(H1)-(H4). Let $\left( X^j,Y^j,Z^j\right) $ ($j=1,2$) be the solutions of
FBSDE (1.1) corresponding to $X_0^1=x_1\in {\bf R}$ and $X_0^2=x_2\in {\bf R}
$, respectively. If $x_1\geq x_2,$ then $X_t^1\geq X_t^2,$ $Y_t^1\geq Y_t^2,$
$\forall t\in \left[ 0,\infty \right] .$
\end{theorem}

\noindent\textbf{Proof.} We set
\begin{eqnarray*}
\tau _x &\doteq &\inf \left\{ t\geq 0;\hat X_t=0\right\} , \\
\tau _y &\doteq &\inf \left\{ t\geq 0;\hat Y_t=0\right\} .
\end{eqnarray*}
Clearly $\tau \leq \tau _x,$ $\tau \leq \tau _y,$ where $\tau $ is defined
in Theorem 2. From Theorem 2, it holds that $\left\langle \hat X_t,\hat
Y_t\right\rangle \geq 0$, $\forall t\geq 0$, and
\begin{eqnarray*}
\hat X_t &=&0,\qquad \forall t\geq \tau _x, \\
\hat Y_t &=&0,\qquad \forall t\geq \tau _y.
\end{eqnarray*}
From the continuity of $\hat X_t$ and $\hat Y_t$, it follows that
$\hat X_t\geq 0$ and $\hat Y_t\geq 0,\quad \forall t\geq 0.$ The
proof is completed.\quad$\Box$

\subsection{Comparison theorem of multi-dimensional BSDEs on infinite horizon
}

We give a Gronwall's lemma which is vital to prove a comparison theorem for
solutions of multi-dimensional infinite horizon BSDEs.

\begin{lemma}
Let $\eta \in L^1\left( \Omega ,{\cal F},P\right) ,\ y_t$ is an adapted
process, and $\varphi \left( s\right) $ is a deterministic function such
that $\int_0^\infty \varphi \left( s\right) ds<\infty $. If $y_t\leq {\bf E}%
\left[ \eta +\int_t^\infty \varphi _sy_sds\left| {\cal F}_t\right. \right] ,$
$\forall t\geq 0,$ then, there exists a positive constant $M$ such that $%
y_t\leq M{\bf E}\left[ \eta \left| {\cal F}_t\right. \right] .$
\end{lemma}

\noindent\textbf{Proof.} From the assumptions, we can get
\begin{equation}
{\bf E}\left[ y_s-\int_s^\infty \varphi _ry_rdr\left| {\cal F}_s\right.
\right] \leq {\bf E}\left[ \eta \left| {\cal F}_s\right. \right] .
\tag{3.2.1}
\end{equation}
Set $\phi _s=\varphi _se^{-\int_s^\infty \varphi _rdr},$ it is easy to check
that $\left( e^{-\int_s^\infty \varphi _rdr}\right) ^{\prime }=\phi _s.$
Multiplying $\phi _s$ on both sides of (3.2.1), it yields
\[
{\bf E}\left[ \phi _sy_s-\phi _s\int_s^\infty \varphi _ry_rdr\left| {\cal F}%
_s\right. \right] \leq {\bf E}\left[ \phi _s\eta \left| {\cal F}_s\right.
\right] .
\]
Integrating on $\left[ t,\infty \right] $, it follows that
\[
\int_t^\infty {\bf E}\left[ \phi _sy_s-\phi _s\int_s^\infty \varphi
_ry_rdr\left| {\cal F}_s\right. \right] ds\leq \int_t^\infty \phi _s{\bf E}%
\left[ \eta \left| {\cal F}_s\right. \right] ds.
\]
Taking conditional expectation under ${\cal F}_t$\ on both sides of the
above inequality
\[
{\bf E}\left[ \int_t^\infty {\bf E}\left[ \phi _sy_s-\phi _s\int_s^\infty
\varphi _ry_rdr\left| {\cal F}_s\right. \right] ds\left| {\cal F}_t\right.
\right] \leq {\bf E}\left[ \int_t^\infty \phi _s{\bf E}\left[ \eta \left|
{\cal F}_s\right. \right] ds\left| {\cal F}_t\right. \right] .
\]
From the property of classical conditional expectation and Fubini theorem,
it follows that
\[
{\bf E}\left[ \int_t^\infty \phi _sy_sds-\int_t^\infty \phi _s\int_s^\infty
\varphi _ry_rdrds\left| {\cal F}_t\right. \right] \leq \int_t^\infty \phi _s%
{\bf E}\left[ \eta \left| {\cal F}_t\right. \right] ds.
\]
After the calculation of integration by parts, it is clear that
\begin{equation}
{\bf E}\left[ \int_t^\infty \varphi _sy_sds\left| {\cal F}_t\right. \right]
\leq \int_t^\infty \varphi _se^{\int_t^s\varphi _rdr}{\bf E}\left[ \eta
\left| {\cal F}_t\right. \right] ds.  \tag{3.2.2}
\end{equation}
Substituting (3.2.2) into the original inequality, the desired
result is obtained.\quad$\Box$

\begin{proposition}
Consider the following backward stochastic differential equations on
infinite horizon,
\begin{equation}
y_t=\xi +\int_t^\infty f\left( s,y_s,z_s\right) ds-\int_t^\infty z_sdW_s,
\tag{3.2.3}
\end{equation}
where $\xi \in L^2,$ $f:\Omega \times {\bf R}^{+}{\bf \times R}^m{\bf \times
R}^{m\times d}{\bf \rightarrow R}^m.$ We assume that

(h.1) for any $\left( y,z\right) \in {\bf R}^m{\bf \times R}^{m\times d},$ $%
f\left( \cdot ,y,z\right) $ is an $\left\{ {\cal F}_t\right\} $%
-progressively measurable process such that ${\bf E}\left( \int_0^\infty
f\left( s,0,0\right) ds\right) ^2<\infty ;$

(h.2) $f$ satisfies Lipschitz condition with Lipschitz function, that is,
there exist two positive deterministic functions $\psi \left( t\right) ,$ $%
\delta \left( t\right) ,$ such that $\left| f\left( t,y^1,z^1\right)
-f\left( t,y^2,z^2\right) \right| \leq \psi \left( t\right) \left|
y^1-y^2\right| +\delta \left( t\right) \left| z^1-z^2\right| ,$ $\forall
\left( t,y^i,z^i\right) \in {\bf R}^{+}{\bf \times R}^m{\bf \times R}%
^{m\times d},$ ($i=1,2$), where $\int_0^\infty \psi \left( s\right)
ds<\infty ,$ $\int_0^\infty \delta ^2\left( s\right) ds<\infty .$ If (h1)
and (h2) hold, then (3.2.3) has a unique solution $\left( y_t,z_t\right) \in
{\cal S}^2\times {\cal H}^2$.
\end{proposition}

The proof is refered to [4].

\bigskip\

Now, we are position to study a comparison theorem for multi-dimensional
BSDEs on infinite horizon.

\begin{theorem}
Consider the following two BSDEs:
\begin{equation}
y_t^1=\xi ^1+\int_t^\infty f^1\left( s,y_s^1,z_s^1\right) ds-\int_t^\infty
z_s^1dW_s  \tag{3.2.4}
\end{equation}
\begin{equation}
y_t^2=\xi ^2+\int_t^\infty f^2\left( s,y_s^2,z_s^2\right) ds-\int_t^\infty
z_s^2dW_s  \tag{3.2.5}
\end{equation}
If $\xi ^i\in L^2,$ $i=1,2,$ and $\xi ^1\geq \xi ^2,$ $f^1,$ $f^2$ enjoy
(h.1) and (h.2), $f^1\left( t,\cdot ,\cdot \right) $ and $f^2\left( t,\cdot
,\cdot \right) $ satisfy the assumption (A2), then $\forall t\geq 0,$ $%
y_t^1\geq y_t^2,$ that is, $y_t^{1,j}\geq y_t^{2,j},$ $1\leq j\leq m.$
\end{theorem}
\noindent\textbf{Proof.}\quad From the above Proposition 5, (3.2.4)
and (3.2.5)
have a unique solution $\left( y^1,z^1\right) $ and $\left( y^2,z^2\right) $%
, respectively. Set
\[
\hat y_t=y_t^2-y_t^1,\text{ }\hat \xi =\xi ^2-\xi ^1,\text{ }\hat f=f^2-f^1,%
\text{ }\hat z=z^2-z^1.
\]
Applying Tanaka-It\^o's formula to $\left| \hat y_t^{+}\right| ^2$ on $%
\left[ t,\infty \right] .$%
\begin{eqnarray*}
&&\left| \left( \hat y_t^j\right) ^{+}\right| ^2+\mathbf{E}^{\mathcal{F}%
_t}\int_t^\infty \sum_{i=1}^d\left| \hat z_s^{j,i}\right|
^2I_{\left\{
y_s^1<y_s^2\right\} }ds \\
&=&2\mathbf{E}^{\mathcal{F}_t}\int_t^\infty \left( \hat y_s^j\right)
^{+}\left[ f^{2,j}\left( s,y_s^1,z_s^2\right) -f^{1,j}\left(
s,y_s^2,z_s^2\right) \right] \mathbf{I}_{\left\{ y_s^1<y_s^2\right\} }ds-%
\mathbf{E}^{\mathcal{F}_t}\int_t^\infty \left( \hat y_s^j\right)
^{+}dL_s^j
\end{eqnarray*}
\begin{equation}
\tag{3.2.6}
\end{equation}

where $L_s^j$ is local time for $y_s^{2,j}-y_s^{1,j}$ at $0$ which claims $%
{\bf E}^{{\cal F}_t}\int_t^\infty \left( \hat y_s^j\right) ^{+}dL_s^j=0.$
Set
\begin{eqnarray*}
\bigtriangleup &=&f^{2,j}\left( s,y_s^2,z_s^2\right) -f^{1,j}\left(
s,y_s^1,z_s^1\right) \\
\ &=&f^{2,j}\left( s,y_s^{2,1},\ldots ,y^{2,j},\ldots
y_s^{2,m},z^{2,1},\ldots ,z^{2,j},\ldots ,z^{2,m}\right) \\
&&\ \ \ \ -f^{1,j}\left( s,y_s^{1,1}+\left( \hat y_s^1\right)
^{+},y_s^{1,2}+\left( \hat y_s^2\right) ^{+},\ldots ,y_s^{2,j},\ldots
,y_s^{1,m}+\left( \hat y_s^m\right) ^{+},z_s^{1,1},\ldots ,z_s^{2,j},\ldots
,z_s^{1,m}\right) \\
&&\ \ \ \ +f^{1,j}\left( s,y_s^{1,1}+\left( \hat y_s^1\right)
^{+},y_s^{1,2}+\left( \hat y_s^2\right) ^{+},\ldots ,y_s^{2,j},\ldots
,y_s^{1,m}+\left( \hat y_s^m\right) ^{+},z_s^{1,1},\ldots ,z_s^{2,j},\ldots
,z_s^{1,m}\right) \\
&&\ \ \ \ -f^{1,j}\left( s,y_s^{1,1},\ldots ,y_s^{1,j},\ldots
,y_s^{1,m},z_s^{1,1},\ldots ,z_s^{1,j},\ldots ,z_s^{1,m}\right) \\
\bigtriangleup _1 &=&f^{2,j}\left( s,y_s^{2,1},\ldots ,y_s^{2,j},\ldots
y_s^{2,m},z_s^{2,1},\ldots ,z_s^{2,j},\ldots ,z_s^{2,m}\right) \\
&&\ \ \ \ -f^{1,j}\left( s,y_s^{1,1}+\left( \hat y_s^1\right)
^{+},y_s^{1,2}+\left( \hat y_s^2\right) ^{+},\ldots ,y_s^{2,j},\ldots
,y_s^{1,m}+\left( \hat y_s^m\right) ^{+},z_s^{1,1},\ldots ,z_s^{2,j},\ldots
,z_s^{1,m}\right) \\
\bigtriangleup _2 &=&f^{1,j}\left( s,y_s^{1,1}+\left( \hat y_s^1\right)
^{+},y_s^{1,2}+\left( \hat y_s^2\right) ^{+},\ldots ,y_s^{2,j},\ldots
,y_s^{1,m}+\left( \hat y_s^m\right) ^{+},z_s^{1,1},\ldots ,z_s^{2,j},\ldots
,z_s^{1,m}\right) \\
&&\ \ \ \ -f^{1,j}\left( s,y_s^{1,1},\ldots ,y_s^{1,j},\ldots
,y_s^{1,m},z_s^{1,1},\ldots z_s^{1,j},\ldots z_s^{1,m}\right) ,
\end{eqnarray*}
so $\bigtriangleup =\bigtriangleup _1+\bigtriangleup _2.$ Obviously,
according to (A2), $\bigtriangleup _1\leq 0.$ From the Lipschitz condition
of $f^1$, we have the following inequality
\[
\bigtriangleup _2\leq \psi \left( s\right) \left( \left| \left( \hat
y_s^1\right) ^{+}\right| +\cdots +\left| \hat y_s^j\right| +\cdots +\left|
\left( \hat y_s^m\right) ^{+}\right| \right) +\delta \left( s\right) \left|
\hat z_s^j\right| .
\]
Immediately, from (3.2.6), it follows that
\begin{eqnarray*}
&&\left| \left( \hat y_t^j\right) ^{+}\right| ^2+\mathbf{E}^{\mathcal{F}%
_t}\int_t^\infty \sum_{i=1}^d\left| \hat z_s^{j,i}\right| ^2\mathbf{I}%
_{\left\{ y_s^1<y_s^2\right\} }ds \\
&\leq &2\mathbf{E}^{\mathcal{F}_t}\int_t^\infty \left[ \left( \hat
y_s^j\right) ^{+}\psi \left( s\right) \left( \left| \left( \hat
y_s^1\right) ^{+}\right| +\cdots +\left| \hat y_s^j\right| +\cdots
+\left| \left( \hat y_s^m\right) ^{+}\right| \right) +\left( \hat
y_s^j\right) ^{+}\delta \left( s\right) \left| \hat z_s^j\right|
I_{\left\{ y_s^1<y_s^2\right\} }\right] ds.
\end{eqnarray*}
Set $\kappa \left( t\right) =\delta ^2\left( t\right) \vee \psi \left(
t\right) ,$ $t\geq 0.$ Obviously, $\int_0^\infty \kappa \left( s\right)
ds<\infty .$ From the inequality $2ab\leq a^2+b^2,$ we get
\begin{eqnarray*}
&&\ \ \ \ \ \left| \left( \hat y_t^j\right) ^{+}\right| ^2+{\bf E}^{{\cal F}%
_t}\int_t^\infty \sum_{i=1}^d\left| \hat z_s^{j,i}\right| ^2{\bf I}_{\left\{
y_s^1<y_s^2\right\} }ds \\
\ &\leq &\left( m+d\right) {\bf E}^{{\cal F}_t}\int_t^\infty \left| \left(
\hat y_s^j\right) ^{+}\right| ^2\kappa \left( s\right) {\bf I}_{\left\{
y_s^1<y_s^2\right\} }ds \\
&&\ \ \ \ \ +{\bf E}^{{\cal F}_t}\int_t^\infty \sum_{i=1}^m\left| \left(
\hat y_s^i\right) ^{+}\right| ^2\kappa (s)I_{\left\{ y_s^1<y_s^2\right\} }ds+%
{\bf E}^{{\cal F}_t}\int_t^\infty \sum_{i=1}^d\left| \hat z_s^{j,i}\right|
^2I_{\left\{ y_s^1<y_s^2\right\} }ds.
\end{eqnarray*}
Finally, there is a sufficiently large constant $K^{\prime }>0,$ such that
\[
\left| \left( \hat y_t^j\right) ^{+}\right| ^2\leq K^{\prime }{\bf E}^{{\cal %
F}_t}\int_t^\infty \sum_{i=1}^m\left| \left( \hat y_s^i\right) ^{+}\right|
^2\kappa \left( s\right) I_{\left\{ y_s^1<y_s^2\right\} }ds,\text{ }1\leq
j\leq m.
\]
Combining all of them, we have
\[
\sum_{j=1}^m\left| \left( \hat y_t^j\right) ^{+}\right| ^2\leq mK^{\prime }%
{\bf E}^{{\cal F}_t}\int_t^\infty \sum_{i=1}^m\left| \left( \hat
y_s^i\right) ^{+}\right| ^2\kappa \left( s\right) {\bf I}_{\left\{
y_s^1<y_s^2\right\} }ds,
\]
by Gronwall's Lemma 4, immediately, we can get
\[
\sum_{j=1}^m\left| \left( \hat y_t^j\right) ^{+}\right| ^2=0,\text{ }\forall
t\in R^{+},
\]
that is $\left( \hat y_t^j\right) ^{+}=0.$ The proof is completed.
\quad$\Box$

\subsection{Multi-dimensional cases of FBSDEs}

From now on, we consider the FBSDE (1.1) whose solution $X$ or $Y$ is
multi-dimensional.

\begin{theorem}
Assume that $n=1,$ $m>1$, and $\left( b,\sigma ,f,\Phi \right) $ satisfy
(H1)-(H4). Let $\left( X^j,Y^j,Z^j\right) $ be the solutions of FBSDE (1.1)
corresponding to $X_0^1=x_1\in {\bf R}$ and $X_0^2=x_2\in {\bf R}$,
respectively. If $x_1\geq x_2,$ $f\left( t,X_t^1,Y,Z\right) $ and $f\left(
t,X_t^2,Y,Z\right) $ satisfy (A2), then $X_t^1\geq X_t^2,$ $Y_t^1\geq Y_t^2,$
$\forall t\in \left[ 0,\infty \right] .$
\end{theorem}

\noindent\textbf{Proof.} Let $\hat X_t=X_t^1-X_t^2,$ $\hat
Y_t=Y_t^1-Y_t^2,$ $\hat Z_t=Z_t^1-Z_t^2,$ we have
\begin{equation}
\left\{
\begin{array}{l}
d\hat X_t=\left( b_t^1\hat X_t+b_t^2\hat Y_t+b_t^3\hat Z_t\right) dt+\left(
\sigma _t^1\hat X_t+\sigma _t^2\hat Y_t+\sigma _t^3\hat Z_t\right) dB_t, \\
-d\hat Y_t=\left( f_t^1\hat X_t+f_t^2\hat Y_t+f_t^3\hat Z_t\right) dt-\hat
Z_tdB_t, \\
\hat X_0=x_1-x_2,\quad \quad \hat Y_\infty =\bar \Phi \cdot \hat X_\infty ,
\end{array}
\right.  \tag{3.3.1}
\end{equation}
where $1\leq j\leq m,$ $l=b,$ $\sigma ,$ $f,$ respectively,
\begin{eqnarray*}
l_t^1 &=&\left\{
\begin{array}{l}
\frac{l\left( t,X_t^1,Y_t^1,Z_t^1\right) -l\left( t,X_t^2,Y_1^1,Z_t^1\right)
}{X_t^1-X_t^2},\text{ \ \ \ if }X_t^1\neq X_t^2, \\
0,\text{\qquad otherwise,}
\end{array}
\right. \\
l_t^2 &=&\left\{
\begin{array}{l}
\frac{l\left( t,X_t^2,Y_t^{21},\ldots ,Y_t^{1j},\ldots
,Y_t^{1m},Z_t^1\right) -l\left( t,X_t^2,Y_t^{21},\ldots ,Y_t^{2j},\ldots
,Y_t^{1m},Z_t^1\right) }{Y_t^{1j}-Y_t^{2j}},\text{\quad if }Y_t^{1j}\neq
Y_t^{2j}, \\
0,\qquad \text{otherwise,}
\end{array}
\right. \\
l_t^3 &=&\left\{
\begin{array}{l}
\frac{l\left( t,X_t^2,Y_t^2,Z_t^{21},\ldots ,Z_t^{1j},\ldots
,Z_t^{1m}\right) -l\left( t,X_t^2,Y_t^2,Z_t^{21},\ldots ,Z_t^{2j},\ldots
,Z_t^{1m}\right) }{Z_t^{1j}-Z_t^{2j}},\text{ \ \ \ if }Z_t^{1j}\neq Z_t^{2j},
\\
0,\text{\qquad otherwise,}
\end{array}
\right. \\
\bar \Phi &=&\left\{
\begin{array}{l}
\frac{\Phi (X_\infty ^1)-\Phi (X_\infty ^2)}{X_\infty ^1-X_\infty ^2},\text{%
\quad if }X_\infty ^1\neq X_\infty ^2, \\
0,\text{\qquad otherwise.}
\end{array}
\right.
\end{eqnarray*}
Now we define a stopping time $\tau \doteq \inf \left\{ t\geq 0;\hat
X_t=0\right\} .$ Consider the following FBSDE:
\begin{equation}
\left\{
\begin{array}{l}
d\hat X_t=\left( b_t^1\hat X_t+b_t^2\hat Y_t+b_t^3\hat Z_t\right) dt+\left(
\sigma _t^1\hat X_t+\sigma _t^2\hat Y_t+\sigma _t^3\hat Z_t\right) dB_t, \\
-d\hat Y_t=\left( f_t^1\hat X_t+f_t^2\hat Y_t+f_t^3\hat Z_t\right) dt-\hat
Z_tdB_t, \\
\hat X_\tau =0.\text{\quad }\hat Y_\infty =\bar \Phi \cdot \hat X_\infty .
\end{array}
\right.  \tag{3.3.2}
\end{equation}
Clearly, (3.3.2) has a unique solution on $\left[ \tau ,\infty
\right] ,$ moreover, $\left( \hat X_t,\hat Y_t,\hat Z_t\right) {\bf
I}_{\left[ \tau ,\infty \right] }\left( t\right) \equiv 0.$ From the
continuity of $\hat X_t$ and $\hat X_0=x_1-x_2\geq 0,$ it follows
that $\hat X_t\geq 0,$ $\forall t\geq 0$. Since $f\left(
t,X_t^1,Y,Z\right) $ and $f\left( t,X_t^2,Y,Z\right) $ satisfy (A2)
and $\bar \Phi \geq 0,$ from Theorem 6, it holds that $Y_t^1\geq
Y_t^2,$ $t\geq 0.$ The proof is completed. \quad$\Box$

Now we study the case for $n>1,$ $m=1.$

\begin{theorem}
Assume that $n>1,$ $m=1$, and $\left( b,\sigma ,f,\Phi \right) $ satisfy
(H1)-(H4). Let $\left( X^j,Y^j,Z^j\right) $ be the solutions of FBSDE (1.1)
corresponding to $X_0^1=x_1\in {\bf R}^n,$ $X_0^2=x_2\in {\bf R}^n$,
respectively. If $x_1\geq x_2,$ $b\left( t,\cdot ,Y^1,Z^1\right) $ and $%
\sigma \left( t,\cdot ,Y^1,Z^1\right) $ satisfy the assumption (A1), then $%
Y_t^1\geq Y_t^2,$ $\forall t\in \left[ 0,\infty \right] .$
\end{theorem}

\noindent\textbf{Proof.}\ Consider the following FBSDE:
\begin{equation}
\left\{
\begin{array}{l}
d\hat X_t=\left( b_t^1\hat X_t+b_t^2\hat Y_t+b_t^3\hat Z_t\right) dt+\left(
\sigma _t^1\hat X_t+\sigma _t^2\hat Y_t+\sigma _t^3\hat Z_t\right) dB_t, \\
-d\hat Y_t=\left( f_t^1\hat X_t+f_t^2\hat Y_t+f_t^3\hat Z_t\right) dt-\hat
Z_tdB_t, \\
\hat X_0=x_1-x_2,\text{\quad }\hat Y_\infty =\bar \Phi \cdot \hat X_\infty .
\end{array}
\right.  \tag{3.3.3}
\end{equation}
where $1\leq j\leq n,$ $l=b,$ $f,$ $\sigma ,$ respectively,
\begin{eqnarray*}
l_t^{1,j} &=&\left\{
\begin{array}{l}
\frac{l\left( t,X_t^{21},\ldots ,X_t^{1j},\ldots
,X_t^{1n},Y_t^1,Z_t^1\right) -l\left( t,X_t^{21},\ldots ,X_t^{2j},\ldots
,X_t^{1n},Y_t^1,Z_t^1\right) }{X_t^{1j}-X_t^{2j}},\text{ \ \ \ if }%
X_{1j}\neq X_{2j}, \\
0,\text{\qquad otherwise,}
\end{array}
\right. \\
l_t^{2,j} &=&\left\{
\begin{array}{l}
\frac{l\left( t,X_t^2,Y_t^1,Z_t^1\right) -l\left( t,X_t^2,Y_t^2,Z_t^1\right)
}{Y_t^1-Y_t^2},\text{ \ \ \ \ if }Y_t^1\neq Y_t^2, \\
0,\text{\qquad otherwise,}
\end{array}
\right. \\
l_t^{3,j} &=&\left\{
\begin{array}{l}
\frac{l\left( t,X_t^2,Y_t^2,Z_t^1\right) -l\left( t,X_t^2,Y_t^2,Z_t^2\right)
}{Z_t^1-Z_t^2},\text{ \ \ \ if }Z_t^1\neq Z_t^2, \\
0,\text{\qquad otherwise,}
\end{array}
\right. \\
\bar \Phi &=&\left\{
\begin{array}{l}
\frac{\Phi \left( X_\infty ^{21},\ldots ,X_\infty ^{1j},\ldots ,X_\infty
^{1n}\right) -\Phi \left( X_\infty ^{21},\ldots ,X_\infty ^{2j},\ldots
,X_\infty ^{1n}\right) }{X_\infty ^{1j}-X_\infty ^{2j}},\text{\quad if }%
X_\infty ^{1j}\neq X_\infty ^{2j}, \\
0,\text{\qquad otherwise.}
\end{array}
\right.
\end{eqnarray*}
Then we use the duality technique to get the dual equations of the above
FBSDE (3.3.3):
\begin{equation}
\left\{
\begin{array}{l}
dP_t=\left( f_t^2P_t-\left( b_t^2\right) ^{*}Q_t-(\sigma
_t^2)^{*}K_t\right) dt+\left( f_t^3P_t-\left( b_t^3\right)
^{*}Q_t-\left( \sigma _t^3\right)
^{*}K_t\right) dB_t, \\
-dQ_t=\left( -\left( f_t^1\right) ^{*}P_t+\left( b_t^1\right)
^{*}Q_t+\left(
\sigma _t^1\right) ^{*}K_t\right) dt-K_tdB_t, \\
P_0=1.\text{\quad }Q_\infty =-\left( \bar \Phi \right) ^{*}P_\infty
.
\end{array}
\right.   \tag{3.3.4}
\end{equation}

It is easy to check that (3.3.4) satisfy (H1)-(H3), so there exists a unique
triple $\left( P_t,Q_t,K_t\right) $ in ${\cal B}^3.$ We define the stopping
time $\tau =\inf \left\{ t\geq 0;\text{ }P_t=0\right\} .$ Consider the
following FBSDE on $\left[ \tau ,\infty \right] ,$%
\begin{equation}
\left\{
\begin{array}{l}
d\bar P_t=\left( f_t^2\bar P_t-\left( b_t^2\right) ^{*}\bar
Q_t-\left( \sigma _t^2\right) ^{*}\bar K_t\right) dt+\left(
f_t^3\bar P_t-\left( b_t^3\right) ^{*}\bar Q_t-\left( \sigma
_t^3\right) ^{*}\bar K_t\right) dB_t,
\\
-d\bar Q_t=\left( -\left( f_t^1\right) ^{*}\bar P_t+\left(
b_t^1\right) ^{*}\bar Q_t+\left( \sigma _t^1\right) ^{*}\bar
K_t\right) dt-\bar K_tdB_t,
\\
\bar P_\tau =0,\text{\quad }\bar Q_\infty =-\left( \bar \Phi \right)
^{*}\bar P_\infty .
\end{array}
\right.   \tag{3.3.5}
\end{equation}
Clearly, (3.3.5) has a unique solution $\left( \bar P_t,\bar Q_t,\bar
K_t\right) $ in ${\cal B}^3.$ Furthermore, $\left( \bar P_t,\bar Q_t,\bar
K_t\right) {\bf I}_{\left[ \tau ,\infty \right] }\left( t\right) \equiv 0.$
Let us set
\begin{eqnarray*}
P_t^{\prime } &=&{\bf I}_{\left[ 0,\tau \right] }\left( t\right) P_t+{\bf I}%
_{\left[ \tau ,\infty \right] }\left( t\right) \bar P_t, \\
Q_t^{\prime } &=&{\bf I}_{\left[ 0,\tau \right] }\left( t\right) Q_t+{\bf I}%
_{\left[ \tau ,\infty \right] }\left( t\right) \bar Q_t, \\
K_t^{\prime } &=&{\bf I}_{\left[ 0,\tau \right] }(t)K_t+{\bf I}_{\left[ \tau
,\infty \right] }(t)\bar K_t.
\end{eqnarray*}
So $\left( P_t^{\prime },Q_t^{\prime },K_t^{\prime }\right) $ is the
solution of (3.3.4). According to $P_0>0$ and the definition of $\tau ,$ we
get $P_t\geq 0.$ Applying It\^o's formula to $\left\langle \hat
X_t,Q_t\right\rangle +\hat Y_tP_{t.}$ on $\left[ 0,\infty \right] ,$ we have
\[
\left\langle \hat X_0,Q_0\right\rangle -\hat Y_0P_0=0.
\]
So $\hat Y_0=-\left\langle \hat X_0,Q_0\right\rangle .$ If we can prove $%
Q_0\leq 0,$ we can obtain $\hat Y_0\geq 0.$ Therefore, we study the
following infinite horizon BSDE
\begin{equation}
\left\{
\begin{array}{l}
-dQ_t=\left( -\left( f_t^1\right) ^{*}P_t+\left( b_t^1\right)
^{*}Q_t+\left(
\sigma _t^1\right) ^{*}K_t\right) dt-K_tdB_t, \\
Q_\infty =-\left( \bar \Phi \right) ^{*}\cdot P_\infty .
\end{array}
\right.   \tag{3.3.6}
\end{equation}

Thanks to (H3) and (H4), we deduce $\left( f_t^1\right) ^{*}>0,$
$\bar \Phi \geq 0.$ By (A1), it holds that
\[
\begin{array}{l}
\frac{b^{1,j}\left( t,x^{2,1},\ldots x^{1,j},\ldots ,x^{1,n},y^1,z^1\right)
-b^{1,j}\left( t,x^{2,1},\ldots ,x^{2,j},\ldots ,x^{1,n},y^1,z^1\right) }{%
x^{1,j}-x^{2,j}}\geq 0, \\
\sigma _{1,t}^{s,j}=0,s\neq j.
\end{array}
\]
We need the following another infinite horizon BSDE
\begin{equation}
\left\{
\begin{array}{l}
-d\tilde Q_t=\left( \left( b_t^1\right) ^{*}\tilde Q_t+\left( \sigma
_t^1\right) ^{*}\tilde K_t\right) dt-\tilde K_tdB_t, \\
\tilde Q_\infty =0.
\end{array}
\right.   \tag{3.3.7}
\end{equation}

Set $f^1\left( P_t,Q_t\right) =-\left( f_t^1\right) ^{*}P_t+\left(
b_t^1\right) ^{*}Q_t+\left( \sigma _t^1\right) ^{*}K_t$ and
$f^2\left( \tilde P_t,\tilde Q_t\right) =\left( b_t^1\right)
^{*}\tilde Q_t+\left( \sigma _t^1\right)
^{*}\tilde K_t.$ Obviously $%
f^1$ and $f^2$ satisfy (A2). From the above comparison Theorem 6, we obtain $%
Q_t\leq 0$, $\forall t\geq 0$. It implies that $\hat Y_0\geq 0,$ immediately.

Set $\tau _y=\inf \left\{ t\geq 0;\text{ }\hat Y_t=0\right\} $. From Theorem
2, we have $\left\langle G\hat X_{\tau _y},\hat Y_{\tau _y}\right\rangle
\equiv 0,$ moreover $\left( \hat X_t,\hat Y_t,\hat Z_t\right) {\bf I}%
_{\left[ \tau _y,\infty \right] }\left( t\right) \equiv 0,$ then we have $%
\hat Y_t\geq 0,$ $\forall t\in \left[ 0,\infty \right] $. The proof
is completed.\quad$\Box$

\bigskip

Lastly, we treat the case that the coefficients of FBSDE $b,$ $f,$ $\Phi ,$
and initial value are all different in two FBSDEs. We can obtain a kind of
comparison theoerm which is only for $Y_0.$

\begin{theorem}
Assume that $n>1,$ $m=1$, and $\left( b^i,\sigma ^i,f^i,\Phi ^i\right) ,$ $%
i=1,$ $2,$ satisfy (H1)-(H4). Let $\left( X^i,Y^i,Z^i\right) $ ($i=1,2$) be
the solutions of FBSDE (1.1) corresponding to $X_0^1=x_1\in {\bf R}^n,$ $%
X_0^2=x_2\in {\bf R}^n,$ respectively. If $x_1\geq x_2,$ $b^1\left(
t,X,Y,Z\right) \geq b^2\left( t,X,Y,Z\right) ,$ $f^1\left( t,X,Y,Z\right)
\geq f^2\left( t,X,Y,Z\right) ,$ $\Phi ^1\left( X\right) \geq \Phi ^2\left(
X\right) ,$ $\sigma ^1=\sigma ^2=\sigma ,$ besides $b^1\left( t,\cdot
,Y^1,Z^1\right) $ and $\sigma ^j\left( t,\cdot ,Y^1,Z^1\right) $, $1\leq
j\leq n,$ satisfy (A1), then $Y_0^1\geq Y_0^2.$
\end{theorem}

\noindent\textbf{Proof.}\ By the similar arguments used in Theorem
8, we give the following FBSDE:
\[
\left\{
\begin{array}{l}
d\hat X_t=\left( b_t^1\hat X_t+b_t^2\hat Y_t+b_t^3\hat Z_t+b^1\left(
t,X_t^2,Y_t^2,Z_t^2\right) -b^2\left( t,X_t^2,Y_t^2,Z_t^2\right) \right) dt
\\
+\left( \sigma _t^1\hat X_t+\sigma _t^2\hat Y_t+\sigma _t^3\hat Z_t\right)
dB_t, \\
-d\hat Y_t=\left( f_t^1\hat X_t+f_t^2\hat Y_t+f_t^3\hat Z_t+f^1\left(
t,X_t^2,Y_t^2,Z_t^2\right) -f^2\left( t,X_t^2,Y_t^2,Z_t^2\right) \right)
dt-\hat Z_tdB_t, \\
\hat X_0=x_1-x_2,\text{\quad }\hat Y_\infty =\hat \Phi \cdot \hat X_\infty
+\Phi ^1\left( X_\infty ^2\right) -\Phi ^2\left( X_\infty ^2\right) ,
\end{array}
\right.
\]
where $l_t^1,$ $l_t^2,$ $l_t^3$ are defined in Theorem 8, and $l=b^1,$ $f^1,$
$\sigma ^1,$ while
\[
\hat \Phi =\left\{
\begin{array}{l}
\frac{\Phi ^1\left( X_t^{21},\ldots ,X_t^{1j},\ldots
,X_t^{1n},Y_t^1,Z_t^1\right) -\Phi ^1\left( t,X_t^{21},\ldots
,X_t^{2j},\ldots ,X_t^{1n},Y_t^1,Z_t^1\right) }{X_{1j}-X_{2j}},\text{ \ \ \
if }X_{1j}\neq X_{2j}, \\
0,\text{\qquad otherwise.}
\end{array}
\right.
\]
Simultaneously, we introduce its dual equation:
\[
\left\{
\begin{array}{l}
dP_t=\left( f_t^2P_t-\left( b_t^2\right) ^{*}Q_t-\left( \sigma
_t^2\right) ^{*}K_t\right) dt+\left( f_t^3P_t-\left( b_t^3\right)
^{*}Q_t-\left( \sigma
_t^3\right) ^{*}K_t\right) dBt, \\
-dQ_t=\left( -\left( f_t^1\right) ^{*}P_t+\left( b_t^1\right)
^{*}Q_t+\left(
\sigma _t^1\right) ^{*}K_t\right) dt-K_tdBt, \\
P_0=1,\text{\quad \quad }Q_\infty =-\left( \hat \Phi \right)
^{*}\cdot P_\infty .
\end{array}
\right.
\]
Applying It\^o's formula to $\left\langle \hat X_t,Q_t\right\rangle +\hat
Y_tP_t$ on $\left[ 0,\infty \right] ,$ we have
\begin{eqnarray*}
\hat Y_0 &=&{\bf E}\left[ \Phi ^1\left( X_\infty ^2\right) -\Phi ^2\left(
X_\infty ^2\right) \right] P_\infty \\
&&\ -{\bf E}\int_0^\infty \left\langle b^1\left( t,X_t^2,Y_t^2,Z_t^2\right)
-b^2(t,X_t^2,Y_t^2,Z_t^2),Q_t\right\rangle dt \\
&&\ +{\bf E}\int_0^\infty \left( f^1\left( t,X_t^2,Y_t^2,Z_t^2\right)
-f^2\left( t,X_t^2,Y_t^2,Z_t^2\right) \right) P_tdt-\left\langle \hat
X_0,Q_0\right\rangle .
\end{eqnarray*}
From $P_t\geq 0$ and $Q_t\leq 0,$ which can be proved similarly in
Theorem 8, associated with $b^1\left( t,X_t^2,Y_t^2,Z_t^2\right)
-b^2\left( t,X_t^2,Y_t^2,Z_t^2\right) \geq 0,$ $f^1\left(
t,X_t^2,Y_t^2,Z_t^2\right) -f^2\left( t,X_t^2,Y_t^2,Z_t^2\right)
\geq 0$. we get that $\hat Y_0\geq 0,$ immediately. The proof is
completed.\quad$\Box$

\end{document}